\documentstyle[12pt]{article}

\newtheorem{prop}{Proposition}

\newcommand {\beq}{\begin{equation}}
\newcommand {\eeq}{\end{equation}}
\newcommand {\beqa}{\begin{eqnarray}}
\newcommand {\eeqa}{\end{eqnarray}}         
\newcommand {\beqs}{\begin{eqnarray*}}
\newcommand {\eeqs}{\end{eqnarray*}}
\newcommand {\bds}{\begin{displaymath}}
\newcommand {\eds}{\end{displaymath}}
\newcommand {\n}{\nonumber\\}
\newcommand {\nn}{\nonumber}

\newcommand{\no}{\noindent}

\newcommand {\bebb}{\begin{thebibliography}}
\newcommand {\eebb}{\end{thebibliography}}      
\newcommand {\bbit}{\bibitem}

\def\al{\alpha}
\def\bt{\beta}
\def\Gm{\Gamma}

\def\dl{\delta}

\def\lm{\lambda}
\def\Lm{\Lambda}

\def\nft{\infty}

\def\rtarr{\rightarrow}

\def\hb2{\frac{\hbar}{2}}
\def\hb{\hbar}

\def\Psd{\Psi ^{\dagger}}
\def\ff(uv){\frac{-(u-v)}{k\hbar}}
\def\f[vu]{\frac{-(v-u)}{k\hbar}}
\def\fw((uv)){\frac{i(u-v)}{k\hbar}}
\def\gw((vu)){\frac{i(v-u)}{k\hbar}}
\def\(#1/#2){\frac{#1}{#2}}
\def\journal#1&#2(#3){\unskip, \sl #1\ \bf #2 \rm(19#3) }
\def\andjournal#1&#2(#3){\sl #1~\bf #2 \rm (19#3) }

\begin{document}


\hfill{AMSS-1999-020}

\baselineskip = 18pt

\vskip 1cm


\begin{center}
{\LARGE\bf Quantum currents in the Coset Space
$SU(2)/U(1)$}

\vspace{1cm}

{\normalsize\bf
Xiang-Mao Ding $^{a,b}${\thanks {E-mail:xmding@itp.ac.cn;
corresponding author}},
Bo-Yu Hou $^c$, Liu Zhao $^c$
}\\
\normalsize $^a$ CCAST, P.O. Box 8730, Beijing,100080, China\\
\normalsize $^b$ Institute of Applied Mathematics,\\
\normalsize Academy of Mathematics and Systems Science,\\
\normalsize Academia Sinica, P.O.Box 2734, 100080, China\\
\normalsize  $^c$ Institute of Modern Physics, Northwest University,
Xian, 710069, China\\
\end{center}

\date{}


\vspace{2cm}

\begin{abstract}

We propose a rational quantum deformed nonlocal currents in
the homogenous space $SU(2)_k/U(1)$, and in terms of it and a 
free boson field a representation for the Drinfeld currents of 
Yangian double at a general level $k=c$ is obtained. In the 
classical limit $\hbar \rightarrow 0$, the quantum nonlocal 
currents become $SU(2)_k$ parafermion, and the realization of 
Yangian double becomes the parafermion realization of 
$SU(2)_k$ current algebra.

\end{abstract}

\vspace{1cm}

\vspace{0.5cm}



\setcounter{section}{1}
\setcounter{equation}{0}
\section*{1.Introduction}

It is well known that Virasoro algebra and affine Lie algebras
play a central role in conformal field theories (CFT) in two
dimensional spacetime\cite{BPZ,Kac}, which in quantum theory corresponds
to fields of massless particles, or massive quantum field theories (QFT)
at the critical points, i.e. when the correlation length becomes infinity
(scale invariant).

It is an interesting problem to describe the symmetry of the QFT off the
critical points. Indeed, impressive progresses have been made in this direction,
part of which is represented by the quantum affine algebra \cite{Drin,JiMi},
which describes the symmetries of the certain lattice models
and, in massive QFT, acts as the role of affine Lie algebra in CFT.
For example, the spin $1/2$ XXZ chain in the region $\Delta<-1$ (off-critical) was
exactly solved with the help of the bosonic representation of
$U_q(\widehat{sl_2})$ at level $1$~\cite{FrJi}; in the same way,
the higher spin $XXZ$ chain cannot be exactly solved without
involving the higher level representation of the quantum affine
algebra $U_q(\widehat{sl_2})$~\cite{Srs}. It is known that a quantum
affine Lie algebra corresponds to certain trigonometric solution of the
quantum Yang-Baxter equation (QYBE); while for a rational solution of
QYBE, the corresponding algebraic structure is the the
quantum double of Yangian with central extension \cite{Drn}.
Contrary to the quantum affine case, the Yangian double with center
is obtained relatively rather late \cite{IK,K,KT} due to some technical
difficulty. Before its explicit formulation was given in \cite{IK,K,KT},
the Yangian double with center was expected to be the symmetry algebra of
quantum non-local currents of massive field theories \cite{Sm}.

Integrable massive QFTs, which
are QFTs away from the critical phase, can be obtained either by perturbations
of CFT~\cite{EgYa} with relevant fields, or in terms of free field realizations,
following Lykyanov~\cite{Lyky}. In fact the free field
realization is a common approach used in both massive quantum field theories
and the representation theory of their dynamical symmetry algebras\cite{Wak}.
The free boson representations of $U_q(\widehat{sl_2})$
with an arbitrary level have been obtained in Refs.\cite{Srs,Mats}.
For the Yangian double with center, the free field representation of
$DY_\hbar(sl_2)$ with level $k(\neq 0, -2)$ was constructed in
\cite{kon}, the level-$1$ and level-$k$ representation of
$DY_{\hbar}(sl_N)$ are obtained in \cite{Ioh} and \cite{DHHZ},
respectively. The level-$k$ free field representation of
$DY_{\hbar}(gl_N)$ is also given in \cite{DHHZ}. However,
all of these are deformations of the Wakimoto modules.
Another kind of module has been known in the representation theories of
the classical \cite{Nem} and quantum \cite{DW,BV}
affine algebras. In the case of Yangian double with center
a free field realization, which corresponds to the Feigin-Fuchs representation
\cite{Nem} of usual affine Lie algebras in the
classical limit, is given in \cite{DZ}. However, different from the corresponding
quantum affine case, the coset structure in that realization is absent.
So, explicit expression of rational quantum currents in the coset space
$SU(2)/U(1)$ is not known yet. In the classical case, the currents
defined on the coset space $SU(2)/U(1)$ subject to a $SU(2)$
nonlocal currents algebra, which is also referred to as parafermion
algebra \cite{Gep,DFSW}. Naively speaking, we would expect
that the rational quantum (or $\hbar$)-deformed
parafermion are relevant to the
off-critical coset WZNW model, or equivalently, the massive
WZNW model in homogeneous space $SU(2)/U(1)$. The $\hbar$-deformed
parafermion is also relevant to the $SU(2)$ Yangian double
with center $DY_{\hbar}(sl_2)_k$, which plays an important role
in condensed state physics \cite{Sm}. In fact the relation
with the $DY_{\hbar}(sl_2)_k$ is our main criterion to
define the rational coset quantum currents.
One of the main properties of
the massive current is nonlocality. It is interesting to study
the rational quantum deformation of coset currents, which is
nonlocal even in the classical case. On the other hand, to
study the rational quantum deformation of the nonlocal algebra
is helpful for the investigating of the quantum
deformation of the chiral vertex operator in more general scheme.

The manuscript is arranged as follows. First, we
briefly review the definition of the $SU(2)$ nonlocal currents in Section 2.
Then we propose a nonlocal quantum currents $\Psi (u)$ and
$\Psi ^{\dagger} (u)$, and the central extension of Yangian
double for the Drinfeld new realization are given in
terms of these rational quantum nonlocal
currents and a $U(1)$ current. In section $3$, we give a bosonic
representation for a nonlocal currents through two sets of bosonic
fields, the operator product expansion (OPE) of this nonlocal
currents satisfies the definition of the proposed quantum nonlocal
currents after a Wick rotation.

\setcounter{section}{2}
\setcounter{equation}{0}
\section*{2.Quantum nonlocal currents}

In this section we propose a quantum nonlocal currents (QNC)
$\Psi (u)$ and $\Psi ^{\dagger} (u) $, in the classical
limit $\hbar \rtarr 0$, this nonlocal quantum currents become the
nonlocal currents $\psi (z)$ and $\psi ^{\dagger} (z)$ of the
coset space $SU(2)/U(1)$ respectively. Then in terms of such
quantum nonlocal currents and a quantum $U(1)$ bosonic current,
we obtain a new kind of representation for the Yangian double
with center in the Drinfeld new realization. Before we going to the details,
we first review briefly the theory of parafermionic
currents in CFT.

Parafermionic currents are primary fields of the 2d CFT. The general
parafermion defined for root lattices are proposed in \cite{DFSW}.
For the case of $SU(2)$, the affine Kac-Moody currents and parafermionic
currents are related by the relations,

\begin{eqnarray}
\chi _{+}(z)&=&\sqrt{k}:\psi (z)\exp(i\phi(z)/\sqrt{k}):,\nonumber\\
\chi _{-}(z)&=&\sqrt{k}:\psi ^{\dagger}(z)\exp(-i\phi(z)/\sqrt{k}):,\nonumber\\
h(z)&=&{i\sqrt{k}}\partial _z \phi (z),
\end{eqnarray}

\noindent where $\chi _{\pm}(z)$ and $h(z)$ are currents of $SU(2)$ affine algebra,
and the radial ordering of the parafermions are given by

\begin{equation}
R\left(\psi _{\alpha}(z)\psi _{\beta}(w)\right)(z-w)^{2\alpha\beta/k}
=R\left(\psi _{\beta}(w)\psi _{\alpha}(z)\right)(w-z)^{2\alpha\beta/k},
\label{eq:rlrr}
\end{equation}

\noindent in which $\al, \, \bt =\pm$. We will drop the $R$ symbol in the
following without of any confusion. The OPE of the parafermionic fields
defined by\cite{DFSW}

\beqa
\psi _{\pm}(z)\psi _{\pm}(w)(z-w)^{2/k}
&=& reg.\nn\\
\psi _{+}(z)\psi _{-}(w)(z-w)^{-2/k}
&=& \frac{1}{(z-w)^2}+reg,
\label{eq:Par}
\eeqa

\no For rational quantum deformation the $SU(2)$ currents algebra
becomes the central extension of Yangian double for the Drinfeld
new realization. From this point view we propose rational
quantum deformation of the above nonlocal currents.
The results are:

\begin{prop}\label{prop1} The rational quantum deformation of nonlocal
currents for $SU(2)/U(1)$ can be defined as:

\beqa
( (u-v)+\hb) \frac{\Gm( \ff(uv)+\(1/k) -1)}{\Gm( \ff(uv)-\(1/k) -1)}
\Psi (u) \Psi(v)
&=&(  (u-v)-\hb)\frac{\Gm( \f[vu] -\(1/k) -1)}{\Gm(\f[vu]-\(1/k) -1)}
\Psi (v) \Psi(u),\nn\\
(  (u-v)-\hb) \frac{\Gm(\ff(uv)+\(1/k))}{\Gm(\ff(uv)-\(1/k) )}
\Psd (u) \Psd (v)
&=&( (u-v)+\hb)\frac{\Gm(\f[vu]+\(1/k) )}{\Gm(\f[vu]-\(1/k) )}
\Psd (v) \Psd (u),\nn\\
\frac{\Gm(\ff(uv)-\(1/k) -\frac{1}{2})}{\Gm(\ff(uv)+\(1/k) -\frac{1}{2})}
\Psi (u) \Psd (v)
&=&\frac{\Gm(\f[vu] -\(1/k) -\frac{1}{2})}{\Gm(\f[vu]+\(1/k) -\frac{1}{2})}
\Psd (v) \Psi(u),
\label{RQP}
\eeqa

\end{prop}

\no At first glance the above relations may seem very strange. Their actual meaning
in the theory of quantum deformed algebras will soon be clear if we resort to
the construction of Yangian double analogous to that of SU(2) current algebra in
terms of ordinary parafermionic currents, and indeed, such a construction exists,
and its classical limit $\hb \rtarr 0$ gives precisely the parafermionic construction
of SU(2) currents.

To make the above statements more concrete, we need to introduce another bosonic
currents.

\beq
[\hat{c}(t),~~\hat{c}(t')]=\frac{\sinh \hb t \sinh
\frac{k}{2} \hb t }{\hb ^2 t}\dl (t+t').
\eeq

\no Defining the bosonic currents

\beqa
C^+(u)=\exp \{-\hb \int ^0 _{-\nft} dt e^{\frac{k}{4} \hb t}
\frac{e^{-iut}}{\sinh \frac{k}{2}\hb t}\hat{c}(t)
-\hb \int ^{\nft} _0 dt e^{-\frac{k}{4} \hb t}
\frac{e^{-iut}}{\sinh \frac{k}{2}\hb t}\hat{c}(t)\},\nn\\
C^-(u)=\exp \{\hb \int ^0 _{-\nft} dt e^{-\frac{k}{4} \hb t}
\frac{e^{-iut}}{\sinh \frac{k}{2}\hb t}\hat{c}(t)
+\hb \int ^{\nft} _0 dt e^{+\frac{k}{4} \hb t}
\frac{e^{-iut}}{\sinh \frac{k}{2}\hb t}\hat{c}(t)\},
\eeqa

\no we have the following exchange relations,

\beqa
\frac{\Gm (\fw((uv)) - \(1/k) -1)}{\Gm (\fw((uv))+ \(1/k) -1)}
C^+(u) C^+(v)&=&C^+(v) C^+(u)
\frac{\Gm (\gw((vu))-\(1/k) -1)}{\Gm (\gw((vu))+\(1/k) -1)},\nn\\
\frac{\Gm(\fw((uv))-\(1/k) )}{\Gm(\fw((uv))+\(1/k) )}C^-(u) C^-(v)&=&C^-(v) C^-(u)
\frac{\Gm(\gw((vu))-\(1/k) )}{\Gm(\gw((vu))+\(1/k) )},\nn\\
\frac{\Gm(\fw((uv))+\(1/k) +\frac{1}{2})}{\Gm(\fw((uv))-\(1/k) +\frac{1}{2})}
C^+(u) C^-(v)&=&C^-(v) C^+(u)
\frac{\Gm(\gw((vu))+\(1/k) +\frac{1}{2})}{\Gm(\gw((vu))-\(1/k) +\frac{1}{2})}.
\eeqa

\no If we perform the Wick rotation for currents $C^{\pm}$, namely
replacing $\hb$ by $-i\hb$, and defining the quantum currents as

\beqa
E(u)=:\Psi (u) C^+(u):,\nn\\
F(u)=:\Psd (u) C^-(u):,
\eeqa

\no we obtain the Drinfeld currents of Yangian double with center, with
relations

\beqa
& & [ H^\pm(u),~H^\pm(v) ] =0,~~~
[c,\mbox{everything}]=0,~~~c=k~,\n
& & (u-v+ \hbar\mp \(k/2) \hbar) (u-v \pm \(k/2) \hbar-\hbar)
H^{\pm}(u) H^{\mp}(v) \n
& &~~~~= (u-v-  \hbar \mp \(k/2) \hbar \hbar)
(u-v \pm \(k/2) \hbar + \hbar)
H^{\mp}(v) H^{\pm}(u),\n
& & (u-v +\hbar \mp \(k/4)\hb) H^{\pm}(u) E(v)
= (u-v - \hbar \mp \(k/4)\hb) E(v) H^{\pm}(u),\n
& & (u-v -\hbar \pm \(k/4) \hbar) H^{\pm}(u) F(v)
= (u-v + \hbar \pm \(k/4)\hbar) F(v) H^{\pm}(u),\n
& & (u-v + \hbar) E(u) E(v)
=  (u-v - \hbar) E(v) E(u), \n
& & (u-v -\hbar) F(u) F(v)
=  (u-v +\hbar) F(v) F(u),\n
& & [ E(u),~F(v) ] =
\frac{1}{\hbar} \left(
\delta( u-v+ \(k/2) \hbar ) H^+(u+\(k/2) \hb ) -
\delta( u-v- \(k/2) \hbar) H^-(v +\(k/2) \hbar) \right),
\label{dyef} \\
\nn
\eeqa

\no wherein

\beqa
H^+(u)=\mbox{exp} \{2\hb \int _0 ^{+\nft} dt~ e^{-iut}\hat{c}(t)\},\nn\\
H^-(u)=\mbox{exp} \{-2\hb \int ^0 _{-\nft} dt~ e^{-iut}\hat{c}(t)\}.
\eeqa

\no Since the fields $C^{\pm}(u)$ and $H^{\pm}(u)$ are only involved
in the quantum $U(1)$ current, the coset structure of classical
currents are preserved, and the coset structure could not be
preserved by another kind of free field representation \cite{DZ}.

\setcounter{section}{3}
\setcounter{equation}{0}
\section*{3.Bosonization of the quantum nonlocal currents}

In order to see that the strange rational deformed quantum
currents are well-defined, i.e. their definition is not
empty, in this section we give a free field realization of
the quantum nonlocal currents. First we introduce two kinds of Heisenberg
algebras,

\beqa
[\hat{b}(t),~~\hat{b}(t')]=-\frac{\sinh \hb t ~ \sinh
\frac{k}{2} \hb t }{\hb ^2 t}\dl (t+t'),\n
~[\hat {\lm}(t),~~\hat{\lm}(t')]=\frac{\sinh \hb t ~ \sinh
\frac{k+2}{2} \hb t }{\hb ^2 t}\dl (t+t').
\eeqa

\no Then defining the following intermediate fields

\beqa
B^+(u)=\mbox{exp} \{-\hb \int ^0 _{-\nft} dt e^{\frac{k}{4} \hb t}
\frac{e^{-iut}}{\sinh \frac{k}{2}\hb t}\hat{b}(t)
-\hb \int ^{\nft} _0 dt e^{-\frac{k}{4} \hb t}
\frac{e^{-iut}}{\sinh \frac{k}{2}\hb t}\hat{b}(t)\},\nn\\
B^-(u)=\mbox{exp} \{\hb \int ^0 _{-\nft} dt e^{-\frac{k}{4} \hb t}
\frac{e^{-iut}}{\sinh \frac{k}{2}\hb t}\hat{b}(t)
+\hb \int ^{\nft} _0 dt e^{+\frac{k}{4} \hb t}
\frac{e^{-iut}}{\sinh \frac{k}{2}\hb t}\hat{b}(t)\},
\eeqa

\beqa
\Lm _+(u)=\mbox{exp} \{-2\hb \int _0 ^{+\nft} dt
\frac{\sinh \frac{\hbar}{2}t }{\sinh \hb t}
e^{-iut}\hat{\lm}(t)\},\nn\\
\Lm _-(u)=\mbox{exp} \{2\hb \int ^0 _{-\nft} dt
\frac{\sinh \frac{\hbar}{2}t }{\sinh \hb t}
e^{-iut}\hat{\lm}(t)\},\nn\\
\bt _+(u)=\mbox{exp} \{-2\hb \int _0 ^{+\nft} dt
\frac{\sinh \frac{\hbar}{2}t }{\sinh \hb t}
e^{-iut}\hat{b}(t)\},\nn\\
\bt _-(u)=exp \{2\hb \int ^0 _{-\nft} dt
\frac{\sinh \frac{\hbar}{2}t }{\sinh \hb t}
e^{-iut}\hat{b}(t)\},
\eeqa

\no we have the following formal commutation relations, i.e. relations to be understood
in the sense of analytic continuation,

\beqa
\bt _+(u)B^{\pm}(v)=&&\frac
{i(u-v)\pm \(k/4) \hb \pm \(1/2) \hb}
{i(u-v)\pm \(k/4) \hb \mp \(1/2) \hb}
B^{\pm}(v)\bt _+(u),\n
B^{\pm}(u) \bt _-(v)=&&\frac
{i(u-v)\pm \(k/4) \hb \pm \(1/2) \hb}
{i(u-v)\pm \(k/4) \hb \mp \(1/2) \hb}
\bt _-(v)B^{\pm}(u),\n
\Lm _+(u)\Lm _-(v)=&&
\frac{\Gm(\frac{i(u-v)}{2\hbar}+\frac{k+2}{4} )
      \Gm(\frac{i(u-v)}{2\hbar}+\frac{k+6}{4} )}
     {\Gm(\frac{i(u-v)}{2\hbar}-\frac{k+2}{4})
     \Gm(\frac{i(u-v)}{2\hbar}-\frac{k-2}{4} )}\n
 &&~~\times
\frac{{\Gm}^2(\frac{i(u-v)}{2\hbar}-\frac{k}{4})}
{{\Gm}^2(\frac{i(u-v)}{2\hbar}+\frac{k+4}{4})}
\Lm _-(v)\Lm _+(u),\n
\bt _+(u)\bt _-(v)=&&
\frac{\Gm(\frac{i(u-v)}{2\hbar}-\frac{k}{4} )
      \Gm(\frac{i(u-v)}{2\hbar}-\frac{k-4}{4} )}
     {\Gm(\frac{i(u-v)}{2\hbar}+\frac{k}{4})
     \Gm(\frac{i(u-v)}{2\hbar}+\frac{k+4}{4} )}\n
 &&~~\times
\frac{{\Gm}^2(\frac{i(u-v)}{2\hbar}+\frac{k+2}{4})}
{{\Gm}^2(\frac{i(u-v)}{2\hbar}-\frac{k-2}{4})}
\bt _-(v)\bt _+(u),
\eeqa

\beqa
\frac{\Gm (\fw((uv)) + \(1/k) -1)}{\Gm (\fw((uv))- \(1/k) -1)}
B^+(u) B^+(v)&=& B^+(v) B^+(u)
\frac{\Gm (\gw((vu))+\(1/k) -1)}{\Gm (\gw((vu))-\(1/k) -1)},\n
\frac{\Gm(\fw((uv))+\(1/k) )}{\Gm(\fw((uv))-\(1/k) )}
B^-(u) B^-(v)&=& B^-(v) B^-(u)
\frac{\Gm(\gw((vu))+\(1/k) )}{\Gm(\gw((vu))-\(1/k) )},\nn\\
\frac{\Gm(\fw((uv))-\(1/k) +\frac{1}{2})}{\Gm(\fw((uv))+\(1/k) +\frac{1}{2})}
B^+(u) B^-(v)&=& B^-(v) B^+(u)
\frac{\Gm(\gw((vu))-\(1/k) +\frac{1}{2})}{\Gm(\gw((vu))+\(1/k) +\frac{1}{2})}.
\eeqa

\no Using these results, we can get a realization of the quantum nonlocal
currents as follows,

\beqa
\Psi (u)=\frac{1}{\hb} :\{\bt _+(u+i\frac{k+2}{4}\hb)\Lm _+(u+i\frac{k}{4}\hb)
-\bt _- (u-i\frac{k+2}{4}\hb)\Lm _-(u-i\frac{k}{4}\hb)\}B^+(u):,\n
\Psd (u)=-\frac{1}{\hb} :\{\bt _+(u-i\frac{k+2}{4}\hb)\Lm _+^{-1}(u-i\frac{k}{4}\hb)
-\bt _-(u+i\frac{k+2}{4}\hb)\Lm _-^{-1}(u+i\frac{k}{4}\hb)\}B^-(u):.
\eeqa

\no By direct calculation, we can show that the OPE of the above nonlocal currents
coincide with the defining relations (\ref{RQP}).

\vskip 1cm

In this paper a set of rational deformed quantum nonlocal $SU(2)$ currents
are proposed, bosonization and their relation with central extension
of Yangian double in the Drinfeld new realization are also discussed.
A related interesting object, i.e. the screening currents for these quantum
nonlocal currents, which is important for the calculation of the
correlation functions, will be considered in a separate paper.

\vskip 1cm

\no {\bf Acknowledgments:}
One of the authors (Ding) would like to thanks Profs. S. K. Wang,
K. Wu and Z. Y. Zhu for fruitful discussion. Ding and Zhao are
supported in part by the "Natural Science Foundation of China".


\bebb{99}

\bibitem{BPZ}
A.A. Belavin. A. M. Polyakov, A. B. Zamolodchikov,
{\it Infinite conformal symmetry in two-dimensional
quantum field theory}
{\em Nucl. Phys.} {\bf B241}, (1984)333-380.

\bibitem{Kac}
V. G. Kac, {\it Infinite Dimensional Lie Algebras}, third ed.,
Cambridge University press, Cambridge 1990.

\bibitem{Drin}
V. G. Drinfeld, In  {\it Proceedings of the International
Congress of Mathematicians}, Berkeley, (1987), 798-.

\bibitem{JiMi}
M. Jimbo and T. Miwa, {\it Algebra analysis of solvable lattice
models}, RIMS {\bf 981},1994.

\bibitem{FrJi}
I.B. Frenkel and N.H. Jing, {\em Proc. Natl. Acad. Sci. USA} {\bf 85},
(1988), 9371.

\bibitem{Srs} 
J. Shiraishi, {\it Free boson representation
$U_q(\hat{sl_2})$ }, {\em Phys. Lett. {\bf A171}} (1992),
243-248.

\bibitem{Drn}
V. G. Drinfeld, {\it Sov. Math. Dokl.} {\bf 36} (1988), 212.

\bibitem{IK}
K. Iohara, M. Kohno, {\it A central extension of
$DY_{\hbar}(gl_2)$ and its vertex representations},
{\em Lett. Math. Phys. {\bf 37}} (1996), 319-328.

\bibitem{K}
S. Khoroshkin, {\it Central Extension of the Yangian Double},
{\em Preprint {\tt q-alg/9602031}}.

\bibitem{KT}
S. Khoroshkin, V. Tolstoy,
{\it Yangian Double}
{\em Lett. Math. Phys. {\bf 36}} (1996), 373-402; hep-th/9406194.

\bibitem{Sm}
F. A. Smirnov, {\it Int. J. Mod. Phys. {\bf A7 suppl. 1B}} (1992),
 813-838, 839.

\bibitem{EgYa}
T.Eguchi and S.-K. Yang, {\it Deformation of conformaal field
theories and soliton equations} {\em Phys. Lett.}
{\bf B 224}, (1989) 373-378; {\it Sine-Gordon theory at rational
values f the coupling constant and minimal conformal models}
{\bf B 235}, (1990), 282-286.

\bibitem{Lyky}
S. L. Lukyanov, Commun. Math. Phys. {\bf 167}, (1995), 183.

\bibitem{Wak}
M. Wakimoto, {\it Commun. Math. Phys. {\bf 104}} (1986), 605.

\bibitem{Mats}
A.Matsuo,  {\it A $q$-deformation of Wakimoto modules,
primary fields and screening operators}, {\em
Commun. Math. Phys.} {\bf 160}, (1994), 33-48.

\bibitem{kon}
H. Konno, {\it Free Field Representation of Level-$k$ Yangian
Double ${\cal D}Y(sl_2)_k$ and Deformation of Wakimoto Modules},
{\em Lett. Math. Phys.} {\bf 40}, (1997), 321-336.

\bibitem{Ioh}
K. Iohara, {\it Bosonic representations of Yangian Double
$DY_{\hbar}(g)$ with $g=gl_N,~sl_N$}
{\em J. Phys.A.: Math. Gen.{\bf 29}},(1997),4593.

\bibitem{DHHZ}
X. M. Ding, B. Y. Hou, B. Yuan Hou, L. Zhao,
{\it Free Boson Representation of \\
$DY_{\hbar}(gl_N)_k$ and  $DY_{\hbar}(sl_N)_k$}
{\em J.Math.Phys.}
{\bf 39}, (1998), 2273-2289.

\bibitem{Nem}
D. Nemeschansky, {\it Feigen-fuchs representation of
$\hat{su(2)_k}$ Kac-Moody algebra}, {\em Phys. Lett.}
{\bf B 224}, (1989)121-124.

\bibitem{DW}
X. M. Ding and P. Wang,
{\it Parafermion representation of the quantum affine
$U_q(\widehat{sl_2})$},
{\em Mod Phys. Lett.} {\bf A11}, (1996), 921.

\bibitem{BV}
Bougourzi, A.H., Vinet, L.
{\it On a bonsic-parafermionic realization of
$U_q(\widehat{sl(2)})$}, {\em Lett. Math. Phys.
{\bf 36}}, (1996),101-108

\bibitem{DZ}
X. M. Ding, and L. Zhao, {\it Free Boson Representation of
$DY_{\hbar}(sl_2)_k$ and the Deformation of the Feigin-Fuchs},
{\em Commun. Theor. Phys.} {\bf 32}, (1999), 103; hep-th/9710106.

\bbit{Gep}
D. Gepner, {\it New conformal fields ftheories associated
with Lie algebras and their partition functions}
{\em Nucl. Phys.} {\bf B290}, (1987), 10-24.

\bibitem{DFSW}
X.M. Ding, H. Fan, K.J. Shi, P. Wang and C.Y. Zhu,
{\it W algebra in SU(3) parafermion model,}
{\em Phys. Rev. Lett.} {\bf 70}, (1993), 2228-2231;
{\it $W_3$ algebra constructed from the SU(3) parafermion,}
{\em Nucl. Phys.} {\bf B422}, (1994), 307.

\eebb

\end{document}